# ROUNDING OF CONTINUOUS RANDOM VARIABLES AND OSCILLATORY ASYMPTOTICS

By Svante Janson

*Uppsala University*

We study the characteristic function and moments of the integer-valued random variable $\lfloor X + \alpha \rfloor$, where $X$ is a continuous random variables. The results can be regarded as exact versions of Sheppard's correction. Rounded variables of this type often occur as subsequence limits of sequences of integer-valued random variables. This leads to oscillatory terms in asymptotics for these variables, something that has often been observed, for example in the analysis of several algorithms. We give some examples, including applications to tries, digital search trees and Patricia tries.

**1. Introduction.** Let $X$ be a continuous random variable with characteristic function $\varphi(t) = \mathbb{E}\, e^{itX}$. We consider the random variable $\lfloor X \rfloor$, that is, $X$ rounded downward to the nearest integer; more generally, we consider $\lfloor X + \alpha \rfloor$ for $\alpha \in \mathbb{R}$. The purpose of this note is to give some formulas for the characteristic function and moments of these rounded variables. (Rounding upward, or to the nearest integer, is a.s. given by replacing $\alpha$ by $\alpha + 1$ and $\alpha + 1/2$, resp., so such roundings are also covered.) One motivation for studying such variables is that they often arise as subsequence limits in distributions of integer-valued random variables; this is discussed in Section 4 and illustrated by several examples from the study of tries, digital search trees and Patricia tries.

Our results can be regarded as exact versions of Sheppard's correction [24]; see Remark 2.5. The results are inspired by some special cases studied in detail using another method by Hitczenko and Louchard [9] and Louchard and Prodinger [17].

**2. Results.** We find it convenient to shift the variables and define, anticipating (4.3) below,

$$(2.1) \qquad X_\alpha := \lfloor X + \alpha \rfloor - \alpha + 1.$$

Alternatively, $X_\alpha = \lceil X + \alpha \rceil - \alpha$ (a.s.) and, letting $\{x\} := x - \lfloor x \rfloor$ denote the fractional part, $X_\alpha = X - \{X + \alpha\} + 1$. Thus $X_\alpha$ is periodic in $\alpha$: $X_{\alpha+n} = X_\alpha$ for $n \in \mathbb{Z}$.

THEOREM 2.1. *Let $X$ be a continuous random variable. Then, with notation as above,*

$$(2.2) \qquad \varphi_{X_\alpha}(t) := \mathbb{E}\, e^{\mathrm{i}t X_\alpha} = \sum_{n=-\infty}^{\infty} e^{2\pi \mathrm{i}n\alpha} \frac{e^{\mathrm{i}t} - 1}{\mathrm{i}(t + 2\pi n)} \varphi(t + 2\pi n),$$

*provided this sum converges; in general,* (2.2) *holds with the sum interpreted as a Cesàro sum. (We interpret the fraction as 1 if $t = -2\pi n$.)*

Proofs are given in Section 3. The sum in (2.2) can be rewritten as

$$(2.3) \qquad \varphi_{X_\alpha}(t) = \sum_{n=-\infty}^{\infty} e^{2\pi \mathrm{i}n\alpha} \widetilde{\varphi}(t + 2\pi n),$$

where

$$\widetilde{\varphi}(t) := \frac{e^{\mathrm{i}t} - 1}{\mathrm{i}t} \varphi(t)$$

is the characteristic function of $X + U$ with $U \sim U(0,1)$ independent of $X$.

REMARK 2.2. Suppose that the moment generating function $\psi(t) := \mathbb{E}\, e^{tX}$ exists in a strip $a < \operatorname{Re} t < b$, with $-\infty \le a < 0 < b \le \infty$. Suppose further that, for example, $\psi(t) = O(|t|^{-\delta})$ for some $\delta > 0$ on each closed substrip $a' \le \operatorname{Re} t \le b'$ with $a < a' < b' < b$. Then $\varphi(t) = \psi(\mathrm{i}t)$, and (2.3) extends by analytic continuation to

$$(2.4) \quad \psi_{X_\alpha}(t) := \mathbb{E}\, e^{t X_\alpha} = \sum_{n=-\infty}^{\infty} e^{2\pi \mathrm{i}n\alpha} \widetilde{\psi}(t + 2\pi n \mathrm{i}), \qquad a < \operatorname{Re} t < b,$$

where $\widetilde{\psi}(t) := \frac{e^t - 1}{t} \psi(t)$ is the moment generating function of $X + U$.

In typical applications, $\varphi$ and its derivatives decrease so rapidly that moments of $X_\alpha$ can be obtained by termwise differentiation in (2.2) or (2.3). We let $D$ denote differentiation. [The $O(|t|^{-\delta})$ condition can be weakened.]

THEOREM 2.3. *If* $\mathbb{E}\,|X|^m < \infty$ *for some* $m \geq 1$ *and* $D^k\varphi(t) = O(|t|^{-\delta})$ *for some* $\delta > 0$ *and* $0 \leq k \leq m-1$, *then*

$$(2.5) \qquad \mathbb{E}\,X_\alpha^m = \mathbb{E}(X+U)^m + \beta_m(\alpha),$$

*where*

$$\beta_m(\alpha) = \sum_{n \neq 0} \mathrm{i}^{-m} D^m \widetilde{\varphi}(2\pi n) \cdot e^{2\pi \mathrm{i} n \alpha}$$

*is a periodic function with mean 0. In particular:*

(i) *If* $\mathbb{E}\,|X| < \infty$ *and* $\varphi(t) = O(|t|^{-\delta})$ *with* $\delta > 0$, *then*

$$(2.6) \qquad \mathbb{E}\,X_\alpha = \mathbb{E}\,X + \frac{1}{2} + \sum_{n \neq 0} \frac{\varphi(2\pi n)}{2\pi n \mathrm{i}} e^{2\pi \mathrm{i} n \alpha}.$$

(ii) *If* $\mathbb{E}\,X^2 < \infty$ *and* $\varphi(t), \varphi'(t) = O(|t|^{-\delta})$ *with* $\delta > 0$, *then*

$$\mathbb{E}\,X_\alpha^2 = \mathbb{E}\left(X + \frac{1}{2}\right)^2 + \frac{1}{12} + \sum_{n \neq 0}\left(\frac{(1 - \mathrm{i}\pi n)\varphi(2\pi n)}{2\pi^2 n^2} - \frac{\varphi'(2\pi n)}{\pi n}\right)e^{2\pi \mathrm{i} n \alpha}$$

*and*

$$(2.7) \qquad \mathrm{Var}\,X_\alpha = \mathrm{Var}\,X + \frac{1}{12} - \sum_{n=1}^{\infty} \frac{|\varphi(2\pi n)|^2}{2\pi^2 n^2} + \tilde{\beta}_2(\alpha),$$

*where* $\tilde{\beta}_2(\alpha) = \beta_2(\alpha) - (2\,\mathbb{E}\,X + 1)\beta_1(\alpha) - \beta_1(\alpha)^2 + \int_0^1 \beta_1(\alpha)^2\,d\alpha$ *is a periodic function with mean 0.*

Similar formulas for higher central moments follow too; see [17] where the cases $m = 2$ and 3 are studied in detail. For comparisons with [9] and [17], note that $X + U - \frac{1}{2}$ has characteristic function $e^{-\mathrm{i}t/2}\widetilde{\varphi}(t) = \frac{\sin(t/2)}{t/2}\varphi(t)$ and moment generating function (if it exists) $\frac{\sinh(t/2)}{t/2}\psi(t)$.

In many cases, $\varphi$ decreases so rapidly that $\beta_m$ is very small; then $X_\alpha$ thus has, for every $\alpha$, approximately the same moments as $X + U$ (for $m$ small, at least); the dependence of $\alpha$ appears only as small oscillations. This is somewhat surprising since each $X_\alpha$ is a discrete variable while $X + U$ is absolutely continuous, but it can be partly explained by the following observation:

REMARK 2.4. If we mix $X_\alpha$ by taking $\alpha$ random and uniform on $[0, 1)$, we obtain exactly the distribution of $X + U$; in other words, $X_U \overset{\mathrm{d}}{=} X + U$. This follows from (2.2), but also (more easily) directly from (2.1) by conditioning on $X$. Hence, if a moment of $X_\alpha$ depends very little on $\alpha$, it has to be close to the corresponding moment of $X + U$.

REMARK 2.5. If we ignore the periodic terms $\beta_m$ in Theorem 2.3, that is, if we approximate the moments of $X_\alpha$ by the moments of $X + U$, we obtain the well-known *Sheppard's corrections* used in statistics for moments of grouped data [24] (see, e.g., [6], Section 27.9). The necessity of the third term in (2.7)—which equals $-\int_0^1 \beta_1(\alpha)^2\,d\alpha = -\operatorname{Var} m(\alpha)$, if $m(\alpha) = \mathbb{E}\,X_\alpha$ and we regard $\alpha$ as random and uniform in $[0,1)$—as an additional correction in order to obtain an unbiased estimate of the variance was pointed out by Carver [4]; see also [5].

REMARK 2.6. The main term $\mathbb{E}(X + U)^m$ in (2.5) is easily expressed in terms of the moments of $X$. One way is through the semi-invariants $\varkappa_m$, for which we have the simple formula

$$\varkappa_m(X + U) = \varkappa_m(X) + \varkappa_m(U) = \varkappa_m(X) + B_m/m, \qquad m \geq 2,$$

where $B_m$ denotes the Bernoulli numbers (cf. [14]).

EXAMPLE 2.7. Let $X = cY$, where $c > 0$ is a constant and $Y$ has the (*Gumbel*) *extreme value distribution* $\mathbb{P}(Y \leq y) = \exp(-\exp(-y))$; thus, $Y = -\ln Z$ with $Z \sim \operatorname{Exp}(1)$. We have

$$(2.8) \qquad \varphi(t) = \mathbb{E}\,e^{\mathrm{i}tcY} = \mathbb{E}\,Z^{-\mathrm{i}tc} = \Gamma(1 - \mathrm{i}ct).$$

Similarly, $\psi(t) = \Gamma(1 - ct) = -ct\Gamma(-ct)$, $\widetilde{\psi}(t) = -c(e^t - 1)\Gamma(-ct)$ and $\widetilde{\varphi}(t) = \widetilde{\psi}(\mathrm{i}t)$. Remark 2.2 applies with $a = -\infty$, $b = 1/c$ and yields

$$\mathbb{E}\,e^{tX_\alpha} = \sum_{n=-\infty}^{\infty} e^{2\pi\mathrm{i}n\alpha}\frac{e^t - 1}{t + 2\pi n\mathrm{i}}\Gamma(1 - c(t + 2\pi n\mathrm{i})), \qquad \operatorname{Re} t < 1/c.$$

For the mean, (2.6) yields, since $\mathbb{E}\,Y = \gamma$,

$$\mathbb{E}\,X_\alpha = c\gamma + \frac{1}{2} - \sum_{n \neq 0} c\Gamma(-2\pi n c\mathrm{i})e^{2\pi\mathrm{i}n\alpha}.$$

Since $|\Gamma(\mathrm{i}t)| = (\pi/(t\sinh(\pi t)))^{1/2} \sim (2\pi)^{1/2}|t|^{-1/2}e^{-\pi|t|/2}$, the terms in the sum are small and decrease rapidly; in the important case $c = 1/\ln 2$, $|\beta_1(\alpha)| < 1.6 \cdot 10^{-6}$ for all $\alpha$. In [17], computations for the first three moments are given, and several other similar examples are also treated. (This function $\beta_1$ occurs in several related contexts too; see, e.g., [8], where $\omega(u) = -\beta_1(u)$, [10], where $p_1(u) = -\beta_1(u)$, [13], Answers 5.2.2-46, 6.3-19 and 6.3-28 ($m = 2$), where $\delta_0(n) = -\beta_1(\log_2 n)$, [18], Theorem 6.2 and Exercise 6.19, where $\delta(u) = -\beta_1(u)$ and [25], page 359, where $P_1(n) = -\beta_1(\log_2 n)$.)

EXAMPLE 2.8. Aldous [1] found in the *random assignment problem* a limit distribution with density $h(x) = e^{-x}(e^{-x} - 1 + x)/(1 - e^{-x})^2$, $x \geq 0$. Let $X$ have this distribution. The mean of the fractional part $\{X\} = X + 1 - X_0$ was used in [19]. The moment generating function of $X$ is

$$\psi(t) = \int_0^\infty \sum_{k=0}^\infty e^{tx-x}(e^{-x} - 1 + x)(k+1)e^{-kx}\,dx$$

$$= \sum_{k=0}^\infty \left( \frac{k+1}{k+2-t} - \frac{k+1}{k+1-t} + \frac{k+1}{(k+1-t)^2} \right)$$

$$= 1 + t\sum_{k=0}^\infty \frac{1}{(k+1-t)^2} = 1 + t\Psi'(1-t), \qquad \operatorname{Re} t < 1,$$

where $\Psi(z) = \Gamma'(z)/\Gamma(z)$ is the logarithmic derivative of the Gamma function; $\Psi'(z) = \sum_{k=0}^\infty (z+k)^{-2}$. It follows that $|\varphi(t)| = |\psi(\mathrm{i}t)| = O(|t|^{-1})$, and (2.6) yields, because $\ln(\Gamma(1+z)\Gamma(1-z)) = \ln(\pi z/\sin(\pi z))$ and thus $\Psi'(1+z) + \Psi'(1-z) = -z^{-2} + \pi^2/\sin^2(\pi z)$,

$$\mathbb{E}\{X\} = \mathbb{E}\,X + 1 - \mathbb{E}\,X_0 = \frac{1}{2} - \sum_{n \neq 0} \frac{1}{2\pi n\mathrm{i}}(1 + 2\pi n\mathrm{i}\Psi'(1 - 2\pi n\mathrm{i}))$$

$$= \frac{1}{2} - \sum_{n=1}^\infty (\Psi'(1 - 2\pi n\mathrm{i}) + \Psi'(1 + 2\pi n\mathrm{i}))$$

$$= \frac{1}{2} - \sum_{n=1}^\infty \left( \frac{1}{4\pi^2 n^2} + \frac{\pi^2}{\sin^2(2\pi^2 n\mathrm{i})} \right)$$

$$= \frac{11}{24} + \sum_{n=1}^\infty \frac{\pi^2}{\sinh^2(2\pi^2 n)},$$

where the last sum is $\approx 2.8 \cdot 10^{-16}$.

In exceptional cases, some oscillating terms may vanish completely.

EXAMPLE 2.9. Let $X \sim U(0, N)$, where $N \geq 1$ is an integer. Then $\varphi(t) = (e^{\mathrm{i}Nt} - 1)/(\mathrm{i}Nt)$. In particular, $\varphi(2\pi n) = 0$ for $n \neq 0$. Hence, Theorem 2.3 yields $\mathbb{E}\,X_\alpha = \mathbb{E}\,X + \frac{1}{2}$, and $\beta_1(\alpha) = 0$. For $0 \leq \alpha \leq 1$, $\lfloor X + \alpha \rfloor \stackrel{\mathrm{d}}{=} \lfloor X \rfloor + Y$, with $Y \sim \mathrm{Be}(\alpha)$ independent of $X$, and a direct calculation easily yields

$$\operatorname{Var} X_\alpha = \operatorname{Var} X - \tfrac{1}{12} + \alpha(1 - \alpha), \qquad 0 \leq \alpha \leq 1.$$

Hence (2.7) holds with $\tilde\beta_2(\alpha) = -\alpha^2 + \alpha - \frac{1}{6}$, $0 \leq \alpha \leq 1$, so $\tilde\beta_2$ does not vanish and is not negligibly small.

EXAMPLE 2.10. Let $N \geq 1$ be an integer and let $X = Y + U = Z + U + U'$, where $Y \sim U(0, N)$, $U, U' \sim U(0, 1)$ and $Z = \lfloor Y \rfloor$ is uniform on $\{0, \ldots, N-1\}$, with $Y, U, U'$ independent. Then $\varphi(t) = -(e^{iNt} - 1)(e^{it} - 1)/(Nt^2)$, and $\varphi(2\pi n) = \varphi'(2\pi n) = 0$ for $n \neq 0$. Consequently, Theorem 2.3 yields $\mathbb{E} X_\alpha = \mathbb{E} X + \frac{1}{2}$ and $\operatorname{Var} X_\alpha = \operatorname{Var} X + \frac{1}{12}$, with $\beta_1(\alpha) = \beta_2(\alpha) = 0$, so there are no oscillations in the first two moments. On the other hand, $\varphi''(2\pi n) = 1/(2\pi n)^2 \neq 0$ and thus $\tilde{\varphi}'''(2\pi n) = 1/(2\pi n)^3 \neq 0$, $n \neq 0$; hence, $\beta_3$ does not vanish and there are oscillations in the third moment $\mathbb{E} X_\alpha^3$.

EXAMPLE 2.11. In Example 4.9, we study a random variable $X$ that has $\varphi(2\pi n) = \psi(2\pi ni) = \Gamma(1 - 2\pi ni/\ln 2)$ exactly as in Example 2.7 with $c = 1/\ln 2$, and thus the oscillating term $\beta_1$ is the same as there.

For the variance, however, there is a surprising cancellation of the oscillations [2, 9, 22]. Indeed, as is shown in Example 4.9, $\operatorname{Var} X_\alpha = 1$ for every $\alpha$, and thus $\tilde{\beta}_2(\alpha)$ in (2.7) vanishes identically; in other words, $\beta_2(\alpha) - (2 \mathbb{E} X + 1)\beta_1(\alpha) - \beta_1(\alpha)^2$ is constant. (Note that $\beta_2$ thus does not vanish.) The $n$th Fourier coefficient of $\tilde{\beta}_2$ is given by Theorem 2.3 and straightforward calculations, using (4.19), (4.20) and (4.21), as

$$\frac{2}{\ln^2 2} \Gamma'\left(-\frac{2\pi n}{\ln 2}i\right) + \frac{2\gamma}{\ln^2 2} \Gamma\left(-\frac{2\pi n}{\ln 2}i\right) - \frac{2}{\ln 2} \sum_{j=1}^{\infty} (-1)^{j-1} \frac{\Gamma(j - 2\pi ni/\ln 2)}{j!(2^j - 1)}$$

$$- \frac{1}{\ln^2 2} \sum_{j \neq 0, n} \Gamma\left(\frac{-2\pi j}{\ln 2}i\right) \Gamma\left(\frac{2\pi(j-n)}{\ln 2}i\right).$$

Prodinger [22] has also given a direct proof by residue calculus of the fact that this rather complicated expression vanishes.

Moreover, (2.7) further shows that

$$\operatorname{Var} X + \frac{1}{12} = 1 + \frac{2}{\ln^2 2} \sum_1^{\infty} \left|\Gamma\left(\frac{2\pi n}{\ln 2}i\right)\right|^2$$

$$= 1 + \frac{1}{\ln 2} \sum_1^{\infty} \frac{1}{n \sinh(2\pi^2 n/\ln 2)}$$

$$\approx 1.000000000001237,$$

another surprise; see [12] and [22]. [$\operatorname{Var} X$ is given by (4.22) below.]

EXAMPLE 2.12. In Example 4.7 below, we find a random variable $X$ with the moment generating function (4.8). In particular, $\varphi(2\pi n) = \psi(2\pi ni) = \Gamma(1 - 2\pi ni/\ln 2)$ is the same as in Example 2.7, and thus the oscillating term $\beta_1$ is the same as there. Since (4.8) implies

$$\mathbb{E} X = \psi'(0) = \eta'(1)/\eta(1) - \Gamma'(1)/\ln 2 = -\alpha + \gamma/\ln 2,$$

with $\alpha := \sum_1^\infty (2^n - 1)^{-1} \approx 1.606695$, Theorem 2.3 yields

$$\mathbb{E}\, X_\alpha = \frac{\gamma}{\ln 2} + \frac{1}{2} - \alpha - \frac{1}{\ln 2} \sum_{n \neq 0} \Gamma\left(-\frac{2\pi n}{\ln 2}\mathrm{i}\right) e^{2\pi \mathrm{i} n\alpha}$$

(cf. the asymptotics for $\mathbb{E}\, U_n$ in [8], [13], Answer 6.3-28, and [18], Theorem 6.2). For higher moments, see [8] and [17], Section 5.5.

Note further that $\eta(2^k) = 0$ by (4.7), and thus (4.8) yields $\psi(k\ln 2 + 2\pi n\mathrm{i}) = 0$ for integers $k \geq 1$ and $n \neq 0$. Hence, (2.4) yields, using (4.9),

$$\begin{align}
(2.9) \qquad \mathbb{E}\, 2^{kX_\alpha} &= \psi_{X_\alpha}(k\ln 2) = \widetilde{\psi}(k\ln 2) = \frac{2^k - 1}{k \ln 2} \psi(k \ln 2) \\
&= \frac{1}{k!} \prod_{j=1}^{k} (2^j - 1)
\end{align}$$

for every integer $k \geq 1$ and every $\alpha$. Thus, as remarked in [17], there is no oscillation in the exponential moments $\mathbb{E}\, 2^{kX_\alpha}$. (There is oscillation for other exponential moments, i.e., for noninteger $k$ and for $k < 0$.)

In other words, the random variables $2^{X_\alpha}$, $0 \leq \alpha < 1$, have the same moments; note that these variables are discrete and supported on disjoint sets. Their mixture $2^{X+U}$ (see Remark 2.4) is a continuous variable with the same moments. We thus have a striking example of distributions *not* determined by their moments.

Example 2.13. In Example 4.8 below, we find a related random variable $X$ with the moment generating function (4.13). In particular, $\varphi(2\pi n) = \psi(2\pi n\mathrm{i}) = \Gamma(1 - 2\pi n\mathrm{i}/\ln 2)/(1 + 2\pi n\mathrm{i}/\ln 2)$, and Theorem 2.3 yields

$$\mathbb{E}\, X_\alpha = \frac{\gamma - 1}{\ln 2} + \frac{3}{2} - \alpha + \frac{1}{\ln 2} \sum_{n \neq 0} \Gamma\left(-1 - \frac{2\pi n}{\ln 2}\mathrm{i}\right) e^{2\pi \mathrm{i} n\alpha}$$

(cf. the asymptotics for $\mathbb{E}\, S_n$ in [16], [13], Answer 6.3-28, and [18], Theorem 6.4). For higher moments, see [17], Section 5.2.

As in Example 2.12, there is no oscillation in the exponential moments $\mathbb{E}\, 2^{kX_\alpha}$; (2.4) and (4.13) yield, in analogy with (2.9),

$$\mathbb{E}\, 2^{kX_\alpha} = \widetilde{\psi}(k\ln 2) = \frac{2^k}{(k+1)!} \prod_{j=1}^{k} (2^j - 1)$$

for every integer $k \geq 1$ and every $\alpha$. Again, the random variables $2^{X_\alpha}$, $0 \leq \alpha < 1$, have different distributions but the same moments.

Remark 2.14. We are studying the distribution of the integer part of $X$ (possibly shifted by a constant). For comparison, note that the fractional part $\{X\}$ has a distribution which is a probability measure on $\mathbb{T} = \mathbb{R}/\mathbb{Z}$ with Fourier coefficients $\varphi(2\pi n)$, numbers that appear frequently in the results above. In particular, we have observed that the $\varphi(2\pi n)$ are the same in Examples 2.7 (with $c = 1/\ln 2$), 2.11 and 2.12, which is equivalent to the fact that the oscillating parts of $\mathbb{E} X_\alpha$ are the same; we now see that this is also equivalent to the fact that $\{X\}$ has the same distribution in these examples.

## 3. Proofs.

PROOF OF THEOREM 2.1. Fix $t$ and consider the periodic function $h(x) := e^{\mathrm{i}t - \mathrm{i}t\{x\}}$ and its Fejér sums $h_N := K_N * h$, where $K_N$ is the Fejér kernel $K_N(x) = \sum_{n=-N}^{N} \widehat{K_N}(n) e^{2\pi\mathrm{i}nx}$ with $\widehat{K_N}(n) = 1 - n/(N+1)$. Note that $|h_N| \le 1$ and that $h_N(x) \to h(x)$ as $N \to \infty$ if $x \notin \mathbb{Z}$ (see, e.g., [26], Theorem III.(3.4)).

We have

$$e^{\mathrm{i}tX_\alpha} = e^{\mathrm{i}t(X - \{X + \alpha\} + 1)} = e^{\mathrm{i}tX} h(X + \alpha),$$

and thus, by dominated convergence,

$$(3.1) \qquad \mathbb{E}\, e^{\mathrm{i}tX_\alpha} = \lim_{N \to \infty} \mathbb{E}(e^{\mathrm{i}tX} h_N(X + \alpha)).$$

The Fourier coefficients of $h$ are

$$\hat{h}(n) = \int_0^1 e^{-\mathrm{i}tx} e^{-2\pi\mathrm{i}nx}\, dx = \frac{e^{\mathrm{i}t} - 1}{\mathrm{i}(t + 2\pi n)}$$

(interpreted as 1 if $t + 2\pi n = 0$). Thus

$$h_N(x) = \sum_{n=-N}^{N} \widehat{K_N}(n) \hat{h}(n) e^{2\pi\mathrm{i}nx} = \sum_{n=-N}^{N} \widehat{K_N}(n) \frac{e^{\mathrm{i}t} - 1}{\mathrm{i}(t + 2\pi n)} e^{2\pi\mathrm{i}nx}$$

and

$$\mathbb{E}(e^{\mathrm{i}tX} h_N(X + \alpha)) = \sum_{n=-N}^{N} \widehat{K_N}(n) \frac{e^{\mathrm{i}t} - 1}{\mathrm{i}(t + 2\pi n)} \mathbb{E}\, e^{\mathrm{i}(t + 2\pi n)X + 2\pi\mathrm{i}n\alpha}$$

$$= \sum_{n=-N}^{N} \widehat{K_N}(n) \frac{e^{\mathrm{i}t} - 1}{\mathrm{i}(t + 2\pi n)} e^{2\pi\mathrm{i}n\alpha} \varphi(t + 2\pi n).$$

The right-hand side is the $N$th Cesàro mean of the sum in (2.2), and the result follows by (3.1) (cf., e.g., [26], Section III.1).  $\square$

PROOF OF THEOREM 2.3. Since $|X_\alpha| \le |X| + 1$, $\mathbb{E}|X_\alpha|^m < \infty$; hence, $\varphi$, $\varphi_{X_\alpha}$ and $\tilde{\varphi}$ are $m$ times continuously differentiable. To see (2.5), we thus have to show that we can differentiate (2.3) termwise $m$ times.

The assumption implies $D^k \widetilde{\varphi}(t) = O(|t|^{-1-\delta})$ for $k \leq m-1$, and thus (2.3) can be differentiated termwise $m-1$ times:

$$(3.2) \qquad D^{m-1} \varphi_{X_\alpha}(t) = \sum_{n=-\infty}^{\infty} e^{2\pi i n \alpha} D^{m-1} \widetilde{\varphi}(t + 2\pi n).$$

By Leibniz's rule,

$$D^{m-1} \widetilde{\varphi}(t) = \frac{e^{it} - 1}{it} D^{m-1} \varphi(t) + \Phi_m(t),$$

where $\Phi_m$ involves $\varphi, \ldots, D^{m-2}\varphi$ and the assumption yields $D\Phi_m(t) = O(|t|^{-1-\delta})$. Hence, for $0 < \varepsilon < 1$ and $n \neq 0$,

$$D^{m-1} \widetilde{\varphi}(2\pi n + \varepsilon) - D^{m-1} \widetilde{\varphi}(2\pi n)$$
$$= \frac{e^{i\varepsilon} - 1}{i(2\pi n + \varepsilon)} D^{m-1} \varphi(2\pi n + \varepsilon) + \Phi_m(2\pi n + \varepsilon) - \Phi_m(2\pi n)$$
$$= O(\varepsilon |n|^{-1-\delta}),$$

so by (3.2) and dominated convergence,

$$i^m \, \mathbb{E} \, X_\alpha^m = D^m \varphi_{X_\alpha}(0) = \lim_{\varepsilon \to 0} \frac{D^{m-1} \varphi_{X_\alpha}(\varepsilon) - D^{m-1} \varphi_{X_\alpha}(0)}{\varepsilon}$$
$$= \sum_{n=-\infty}^{\infty} e^{2\pi i n \alpha} D^m \widetilde{\varphi}(2\pi n).$$

The term with $n = 0$ is $D^m \widetilde{\varphi}(0) = i^m \, \mathbb{E}(X + U)^m$, and (2.5) follows.

For $n \neq 0$, $\widetilde{\varphi}(2\pi n) = 0$, $\widetilde{\varphi}'(2\pi n) = \varphi(2\pi n)/(2\pi n)$ and

$$\widetilde{\varphi}''(2\pi n) = i \frac{\varphi(2\pi n)}{2\pi n} - 2 \frac{\varphi(2\pi n)}{(2\pi n)^2} + 2 \frac{\varphi'(2\pi n)}{2\pi n},$$

and the formulas for $\beta_1$ and $\beta_2$ used in (i) and (ii) follow; recall also that $\operatorname{Var} U = \frac{1}{12}$.

Finally, $\int_0^1 |\beta_1(\alpha)|^2 \, d\alpha$ is evaluated by Parseval's relation. $\quad\square$

## 4. Asymptotics of integer-valued random variables.

A common situation in the study of integer-valued random variables is the following: $Y_1, Y_2, \ldots$ is a sequence of integer-valued random variables, $a_1, a_2, \ldots$ is a sequence of real numbers and $X$ is a random variable such that

$$(4.1) \qquad \mathbb{P}(Y_n - a_n \leq x_n) = \mathbb{P}(X \leq x_n) + o(1) \qquad \text{as } n \to \infty,$$

for every sequence $x_n$ such that $x_n + a_n \in \mathbb{Z}$; equivalently, for every sequence $k_n$ of integers,

$$(4.2) \qquad \mathbb{P}(Y_n \leq k_n) = \mathbb{P}(X + a_n \leq k_n) + o(1) \qquad \text{as } n \to \infty.$$

Some examples are given below and in, for example, [17] and [18]; in particular, [18], Figures 5.5 and 6.3, provide nice illustrations.

Note that if (4.1) would hold for any real numbers $x_n$, we would have $Y_n - a_n \xrightarrow{\mathrm{d}} X$; however, for integer-valued variables, such convergence in distribution is possible only if the fractional parts $\{a_n\}$ converge as elements of the circle $\mathbb{R}/\mathbb{Z}$ (i.e., we identify 0 and 1). In general, it is thus impossible to have convergence in distribution, and (4.1) is a natural substitute. For example, (4.1) should be expected for a sequence $Y_n$ of random variables that arise in some natural way and have bounded variances. [Of course, if the variance tends to $\infty$, we can have convergence of $(Y_n - a_n)/\operatorname{Var}(Y_n)^{1/2}$ to a continuous random variable; we do not consider that case in this paper.] Note that the centering constants $a_n$ can be regarded as approximations of the median (and typically of the mean as well).

Before proceeding, let us note that (4.1) and (4.2) can appear in different forms. For convenience, we state some of these versions in the following simple lemma; the proof is left to the reader. [$d_{\mathrm{TV}}$ denotes the total variation distance; for integer-valued random variables $X$ and $Y$, $d_{\mathrm{TV}}(X,Y) = \frac{1}{2}\sum_k |\mathbb{P}(X=k) - \mathbb{P}(Y=k)|$.]

LEMMA 4.1.   *If $Y_n$ are integer-valued random variables, $a_n$ are real numbers and $X$ is a random variable, the following are equivalent:*

  (i)  *(4.1) holds for every sequence $x_n$ such that $x_n + a_n \in \mathbb{Z}$.*
  (ii)  *(4.1) holds for every bounded sequence $x_n$ such that $x_n + a_n \in \mathbb{Z}$.*
  (iii)  *(4.2) holds for every sequence $k_n$ of integers.*
  (iv)  *(4.2) holds for every sequence $k_n$ of integers such that $k_n = a_n + O(1)$.*
  (v)  $\sup_{x \in \mathbb{Z} - a_n} |\mathbb{P}(Y_n - a_n \le x) - \mathbb{P}(X \le x)| \to 0$ *as $n \to \infty$.*
  (vi)  $\sup_{k \in \mathbb{Z}} |\mathbb{P}(Y_n \le k) - \mathbb{P}(X + a_n \le k)| \to 0$ *as $n \to \infty$.*
  (vii)  $d_{\mathrm{TV}}(Y_n, \lceil X + a_n \rceil) \to 0$ *as $n \to \infty$.*

Although (4.1) typically prohibits convergence in distribution, it implies convergence in distribution of suitable subsequences. Indeed, suppose that $X$ is a continuous random variable. It is easily checked that if (4.1) holds, then, along any subsequence such that $\{a_n\} \to \alpha \in [0,1]$,

$$(4.3) \qquad\qquad Y_n - a_n \xrightarrow{\mathrm{d}} \lceil X + \alpha \rceil - \alpha = X_\alpha.$$

Conversely, if (4.3) holds for every subsequence with $\{a_n\} \to \alpha$, then (4.1) holds.

REMARK 4.2.   Lemma 4.1(vii) can be written $d_{\mathrm{TV}}(Y_n - a_n, X_{\{a_n\}}) \to 0$. If, further, $m \ge 1$ and $(Y_n - a_n)^m$ is uniformly integrable, this implies $\mathbb{E}(Y_n - a_n)^m = \mathbb{E}\, X_{\{a_n\}}^m + o(1)$, so asymptotics for moments of $Y_n$ follow from Theorem 2.3. Uniform integrability is easily verified in many examples, but for simplicity we omit discussions of it in our examples.

Example 4.3. A simple example is the *maximum of i.i.d. geometric random variables*. For convenience, we shift the variables by 1; thus, let $Y_n = \max(Z_1, \ldots, Z_n)$, where $Z_i$ are independent and $Z_i - 1 \sim \mathrm{Ge}(p)$ with $0 < p < 1$, that is, $\mathbb{P}(Z_i = k) = p(1-p)^{k-1}$, $k = 1, 2, \ldots$. With $q = 1 - p$ and $a_n = \log_{1/q} n$, it is easily seen that if $x_n + a_n \in \mathbb{Z}_+$, then

$$\mathbb{P}(Y_n \le x_n + a_n) = (1 - q^{x_n + a_n})^n$$
$$= (1 - q^{x_n}/n)^n$$
$$= \exp(-q^{x_n}) + o(1),$$

uniformly in such $x_n$, and (4.1) holds with $X$ having the extreme value distribution $\mathbb{P}(X \le x) = \exp(-q^x)$; this is the distribution in Example 2.7 with $c = 1/|\ln q|$. Hence, (4.3) also holds with this $X$. See further [15], Chapter 1.

Example 4.4. Let $D_n$ be the *depth of a trie* [25], Section 1.1, constructed from $n$ independent random infinite strings of the $m$ symbols $0, \ldots, m-1$, where the symbols in each string are chosen independently and uniformly. It is easily shown (see, e.g., [20] or [18], Section 5.5) that

$$\mathbb{P}(D_n \le k) = (1 - m^{-k})^{n-1}, \qquad k \ge 1.$$

$D_n$ thus have the same distribution as $Y_{n-1}$ in Example 4.3, with $q = 1/m$. Hence, (4.1) and (4.3) hold with $a_n = \log_m n$ and $X$ having the extreme value distribution in Example 2.7 with $c = 1/\ln m$. See further [18], Section 5.5 and [17], Sections 4.1 and 5.1.

Example 4.5. Let $H_n$ be the *height (i.e., maximum depth) of the trie* in Example 4.4. Then (see, e.g., [7, 20], [3], Example 6.2.1, [18], Section 5.6) (4.1) and (4.3) hold with $a_n = 2\log_m n - \log_m 2$ and $X$ as in Example 4.4 having the extreme value distribution in Example 2.7 with $c = 1/\ln m$.

Several other similar (but sometimes more complicated) examples from analyses of algorithms are studied in, for example, [17, 18, 25] and the references given therein.

These references also contain related examples where it is not obvious whether they fit in the above framework or not. Note that if $F$ is the distribution function of $X$, then (4.1) and (4.2) are equivalent to

$$(4.4) \qquad \mathbb{P}(Y_n \le k_n) = F(k_n - a_n) + o(1) \qquad \text{as } n \to \infty,$$

for every sequence $k_n$ of integers. However, there are several cases where (4.4) has been proved for some continuous function $F : \mathbb{R} \to [0, 1]$ that satisfies $\lim_{x \to -\infty} F(x) = 0$ and $\lim_{x \to \infty} F(x) = 1$, but such that it is not obvious that $F$ is monotone. [If it is, it is a distribution function, and thus (4.1) holds for $X$ having this distribution.]

In some cases, this problem is solved by the following lemma:

LEMMA 4.6. *Let $Y_n$ be a sequence of integer-valued random variables that is monotone: $Y_n \leq Y_{n+1}$ for $n \geq 1$. Suppose further that* (4.4) *holds for a right-continuous function $F$ with $\lim_{x \to -\infty} F(x) = 0$ and $\lim_{x \to \infty} F(x) = 1$, and a sequence $a_n$ such that $a_n \to +\infty$ and $a_{n+1} - a_n \to 0$ as $n \to \infty$. Then $F$ is a distribution function.*

PROOF. We have to show that $F$ is monotone. Let $x < y$ be two real numbers. For an integer $k > a_0 + y$, let $n_k := \max\{n : a_n < k - x\}$ and $n'_k := \max\{n : a_n < k - y\}$. The assumptions imply that $n_k \geq n'_k$ and that, as $k \to \infty$, $n_k, n'_k \to \infty$ and $0 < (k - x) - a_{n_k} \leq a_{n_k+1} - a_{n_k} \to 0$ and, similarly, $0 < (k - y) - a_{n'_k} \to 0$. Hence, it follows from (4.4) that, as $k \to \infty$,

$$
\begin{aligned}
(4.5) \qquad \mathbb{P}(Y_{n_k} \leq k) &= F(k - a_{n_k}) + o(1) \to F(x), \\
\mathbb{P}(Y_{n'_k} \leq k) &= F(k - a_{n'_k}) + o(1) \to F(y).
\end{aligned}
$$

Furthermore, $Y_{n_k} \geq Y_{n'_k}$, and thus $\mathbb{P}(Y_{n_k} \leq k) \leq \mathbb{P}(Y_{n'_k} \leq k)$ for each $k$. Consequently, (4.5) implies $F(x) \leq F(y)$. Thus $F$ is monotone. □

EXAMPLE 4.7. Let $U_n$ be a pure birth process with $U_0 = 0$ and $\mathbb{P}(U_{n+1} = k + 1 \mid U_n = k) = 2^{-k}$. This process has been studied by [8] as *approximate counting*. The same random variables (more precisely, random variables with the same distributions) have also been studied as the number of comparisons in an *unsuccessful search in a random binary digital search tree* with $n$ records; see [13, 18] or [25] for definitions. (To see the equivalence, construct a digital search tree from $n$ random binary strings, with bits i.i.d. and Be(1/2), and search (a.s. unsuccessfully) for another random binary string $w$. Say that this requires $k$ comparisons. If the tree is enlarged by a new string $w_{n+1}$, the search depth of $w$ increases from $k$ to $k + 1$ if and only if the first $k$ bits of $w$ and $w_{n+1}$ coincide. See also [21].)

It is shown in [8], Proposition 3 (in an equivalent form) and in [18], Section 6.4, that (4.4) holds for $U_n$ with $a_n = \log_2 n$ and

$$
(4.6) \qquad F(x) = \sum_{j=0}^{\infty} \frac{R_j}{Q} \exp(-2^{-(x-j)}),
$$

where $R_j := (-1)^j \prod_{k=1}^{j} (2^k - 1)^{-1}$ and $Q = \sum_{j=0}^{\infty} R_j$. Thus $F$ is an infinite linear combination of translates of the extreme value distribution function of $Y/\ln 2$ in Example 2.7; note, however, that some $R_j$ are negative, so this is not a simple mixture of these distributions. This is a case where Lemma 4.6 applies, and thus $F$ is a distribution function. [It seems difficult to give an analytic proof that the sum in (4.6) is monotone.] Let $X$ have this distribution.

Using the identity [16, 18]

$$(4.7) \qquad \sum_{j=0}^{\infty} R_j z^j = \eta(z) := \prod_{k=1}^{\infty} \left( 1 - \frac{z}{2^k} \right), \qquad z \in \mathbb{C},$$

it follows from (2.8) that $\varphi(t) = Q^{-1} \eta(e^{it}) \Gamma(1 - it/\ln 2)$. Moreover, $Y$ has the moment generating function $\psi_Y(t) = \Gamma(1 - t/\ln 2)$ in the region $\operatorname{Re} t < \ln 2$, and it follows that $X$ has, in this domain, the moment generating function

$$(4.8) \qquad \psi(t) = \sum_{j=0}^{\infty} \frac{R_j}{Q} e^{jt} \Gamma(1 - t/\ln 2) = \frac{\eta(e^t)}{\eta(1)} \Gamma(1 - t/\ln 2).$$

The function $\eta$ is entire, with zeros at $2^k$, $k = 1, 2, \dots$, and thus $\eta(e^t)$ has a zero at every pole of $\Gamma(1 - t/\ln 2)$, since these are at $t = k \ln 2$, $k = 1, 2, \dots$. Consequently, the right-hand side of (4.8) is entire, which implies that the moment generating function exists for all $t$ and is given by (4.8). [This shows that there is significant cancellation in (4.6); the tail $1 - F(x)$ decreases super-exponentially as $x \to \infty$, while the individual terms in the sum approach their limits $R_j/Q$ at the rate $2^{-x}$.] The values of $\psi(t)$ at the poles of $\Gamma(1 - t/\ln 2)$ are easily calculated: For $k = 1, 2, \dots$,

$$(4.9) \quad \psi(k \ln 2) = \frac{\eta'(2^k)}{\eta(1)} 2^k \left( -\frac{1}{\ln 2} \right)^{-1} \operatorname{Res}_{1-k} \Gamma = \frac{\ln 2}{(k-1)!} \prod_{j=1}^{k-1} (2^j - 1).$$

It follows by the results of Example 2.12 and (easily verified) uniform integrability that there are oscillatory terms in the asymptotics for $\mathbb{E} U_n$—see [8], [13], Answer 6.3-28, and [18], Theorem 6.2—but not for $\mathbb{E} 2^{kU_n}$ when $k$ is a positive integer. Indeed, it is easy to find exact expressions for these exponential moments by recursion (at least for small $k$); they are polynomials in $n$, and the first three are $\mathbb{E} 2^{U_n} = n + 1$, $\mathbb{E} 2^{2U_n} = 3\binom{n+1}{2} + 1$, $\mathbb{E} 2^{3U_n} = 21\binom{n+1}{3} + 7n + 1$ (cf. [8], Proposition 0).

Note that the moments of $2^{U_n}$ thus converge to the values in (2.9), although the distribution does not converge. Thus the method of moment fails for these variables.

Finally, we remark that the problems can be generalized to other bases by taking instead $\mathbb{P}(U_{n+1} = k + 1 \mid U_n = k) = a^{-k}$, where $a > 1$ is an arbitrary real for approximate counting (see [8]) and $a = m$ is an integer for $m$-ary digital search trees; see [13]. Similar results hold for the generalization, replacing 2 by $a$ in (4.6), (4.7), (4.8), (4.9) and elsewhere.

EXAMPLE 4.8. Let $S_n$ be the number of comparisons in a *successful search in a random binary digital search tree* with $n$ records [13, 16, 18]. Then,

$$(4.10) \qquad S_n \overset{\mathrm{d}}{=} 1 + U_{I_n},$$

where $(U_n)_0^\infty$ is as in Example 4.7, and $I_n$ is uniform on $\{0, \ldots, n-1\}$ and independent of $(U_n)$. Louchard [16] has proved (4.4) with the limiting function

$$(4.11) \quad F(x) = 2^x \left(1 - \frac{1}{2Q} \sum_{i=0}^\infty \frac{R_i}{2^i} e^{-2^{-x+i+1}}\right) = \frac{1}{Q} \sum_{i=0}^\infty \frac{R_i}{2^{i+1-x}} (1 - e^{-2^{-x+i+1}}).$$

Since (4.10) implies that $S_n$ is (stochastically) monotone, Lemma 4.6 applies, and this function is a distribution function.

Alternatively, we may use (4.10) and the results of Example 4.7. If we let $X_u$ denote the random variable with distribution function given by (4.6), and still let $a_n = \log_2 n$, then (4.10) and (4.2) yield, for every sequence $k_n = a_n + x_n$ of integers, as $n \to \infty$,

$$(4.12) \quad \begin{aligned} \mathbb{P}(S_n \le k_n) &= \frac{1}{n} \sum_{m=0}^{n-1} \mathbb{P}(U_m \le k_n - 1) \\ &= \frac{1}{n} \sum_{m=0}^{n-1} \mathbb{P}(X_u + a_m \le k_n - 1) + o(1) \\ &= \frac{1}{n} \sum_{m=0}^{n-1} \mathbb{P}(X_u + a_m - a_n + 1 \le x_n) + o(1) \\ &= \mathbb{P}(X_u + a_{I_n} - a_n + 1 \le x_n) + o(1), \end{aligned}$$

with $I_n$ independent of $X_u$. Since $I_n/n \xrightarrow{d} U \sim U(0,1)$, we have

$$a_{I_n} - a_n = \log_2(I_n/n) \xrightarrow{d} \log_2 U = -Z/\ln 2,$$

with $Z := -\ln U \sim \text{Exp}(1)$. Hence, (4.12) yields (4.1) with $X := X_u - Z/\ln 2 + 1$, where $Z$ and $X_u$ are independent. [This $X$ thus must have the distribution function (4.11), which also can easily be verified from (4.6).] By (4.8), $X$ has the moment generating function

$$(4.13) \quad \psi(t) = \psi_{X_u}(t) \frac{e^t}{1 + t/\ln 2} = e^t \frac{\eta(e^t)\Gamma(1 - t/\ln 2)}{\eta(1)(1 + t/\ln 2)}, \qquad \text{Re}\, t > -\ln 2.$$

In analogy with the results above for unsuccessful search, there is oscillation in the asymptotics for $\mathbb{E}\, S_n$—see [16], [13], Answer 6.3-28, [18], Theorem 6.4, and [17], Section 5.2—but not for the exponential moments $\mathbb{E}\, 2^{kS_n}$, which are polynomials in $n$ by (4.10) and the corresponding result for $U_n$.

EXAMPLE 4.9. Let, for $s > 0$, $Y_s := \sum_{k=1}^\infty I_{sk}$, where the $I_{sk} \sim \text{Be}(1 - e^{-s/2^k})$ are independent. Thus $Y_s$ has the probability generating function

$$D(s, u) := \mathbb{E}\, u^{Y_s} = \prod_{k=1}^\infty (e^{-s/2^k} + u(1 - e^{-s/2^k})).$$

These variables occur in the study of the *depth in a random binary Patricia trie*; see [13] or [25] for definitions. Indeed, if $D_n$ is this depth in a Patricia trie with $n$ records, then $D_n$ can be obtained by taking a random binary string $w$ and comparing it to $n-1$ other random binary strings; $D_n$ equals the number of indices $k$ such that there is at least one of the strings that is equal to $w$ in the first $k-1$ positions but differs in the $k$th. If we Poissonize and instead compare with $\mathrm{Po}(s)$ other strings, we obtain $Y_s$. [The number of strings that first differ in the $k$th position is $\mathrm{Po}(s/2^k)$, and these numbers are independent.] Note further that depoissonization is easy because the variables are monotone: if $n_\pm = n \pm n^{2/3}$ and $N_{n\pm} \sim \mathrm{Po}(n_\pm)$, then $\mathbb{P}(N_{n-} \le n-1 \le N_{n+}) \to 1$ and thus $\mathbb{P}(Y_{n-} \le D_n \le Y_{n+}) \to 1$; hence, the asymptotics for $Y_s$ below yield the same results for $D_n$.

The same variables occur as the *number of different values in a sample of $n$ values from a geometric* $\mathrm{Ge}(1/2)$ *distribution* [2] (this is our $D_{n+1}$) and in the study of *distinct parts in a random decomposition* [9].

Let $B_m(s) := \mathbb{P}(Y_s = m)$; thus $0 \le B_m(s) \le 1$ and $D(s,u) = \sum_m B_m(s)u^m$. It is shown in [23] by a manipulation of generating functions that (see [23], (35))

$$\mathbb{P}(Y_s \le k) = \sum_{m=0}^{k} B_m(s/2^{k+1-m})e^{s/2^{k+1-m}},$$

and thus, for any bounded $x_s$ with $\log_2 s + x_s \in \mathbb{Z}$,

$$(4.14) \qquad \mathbb{P}(Y_s - \log_2 s \le x_s) = F(x_s) + o(1) \qquad \text{as } s \to \infty,$$

with

$$(4.15) \qquad F(x) = \sum_{m=0}^{\infty} B_m(2^{m-1-x})e^{-2^{m-1-x}}, \qquad -\infty < x < \infty.$$

Thus (4.4) holds with $a_n = \log_2 n$.

It is easily seen that each $B_m(s)$ is a continuous function, and thus $F(x)$ is continuous (by dominated convergence) with $F(x) \to 0$ as $x \to -\infty$ (by dominated convergence) and $F(x) \to 1$ as $x \to \infty$ [because $F(x) \le 1$ by (4.14) and $F(x) \ge B_0(2^{-1-x})e^{-2^{-1-x}} = e^{-2^{-x}}$]. Lemma 4.6 shows that $F$ is a distribution function of some random variable $X$.

To see that $X$ has an entire moment generating function, note first that for $t > 0$,

$$\mathbb{E}\,e^{t(Y_s - \log_2 s)} \le \mathbb{E} \prod_{k > \log_2 s} e^{tI_{sk}} = \prod_{k > \log_2 s} \left(1 + (e^t - 1)(1 - e^{-s/2^k})\right)$$

$$\le \exp\left((e^t - 1) \sum_{k > \log_2 s} s/2^k\right) \le \exp(2(e^t - 1)),$$

so by (4.3) and Fatou's lemma, $\mathbb{E}\, e^{tX_\alpha} < \infty$ and thus $\mathbb{E}\, e^{tX} < \infty$. The case $t < 0$ can be treated similarly, but we will instead perform an exact computation. If $u > 0$, then, by Fubini and (4.15), with $x = y + m - 1$ and $s = 2^{-y}$,

$$
\begin{aligned}
\psi(-u) := \mathbb{E}\, e^{-uX} &= \mathbb{E}\, u \int_X^\infty e^{-ux}\, dx = u \int_{-\infty}^\infty e^{-ux} F(x)\, dx \\
&= u \sum_{m=0}^\infty \int_{-\infty}^\infty e^{-ux} B_m(2^{m-1-x}) e^{-2^{m-1-x}}\, dx \\
&= u \sum_{m=0}^\infty e^{-u(m-1)} \int_{-\infty}^\infty e^{-uy} B_m(2^{-y}) e^{-2^{-y}}\, dy \\
&= u e^u \int_{-\infty}^\infty e^{-uy} D(2^{-y}, e^{-u}) e^{-2^{-y}}\, dy \\
&= \frac{u e^u}{\ln 2} \int_0^\infty s^{u/\ln 2 - 1} D(s, e^{-u}) e^{-s}\, ds < \infty,
\end{aligned}
$$

since $0 \le D(s, e^{-u}) \le 1$. Hence, $\psi$ is an entire function with, by analytic continuation,

$$
(4.16) \qquad \psi(t) = -\frac{t e^{-t}}{\ln 2} \int_0^\infty D(s, e^t) s^{-t/\ln 2 - 1} e^{-s}\, ds, \qquad \operatorname{Re} t < 0.
$$

Using integration by parts, we can obtain similar formulas (involving also $\partial D / \partial s$ and possibly higher derivatives) that extend into the right half-plane; in particular, this yields formulas for the characteristic function $\varphi(t) = \psi(\mathrm{i}t)$ that can be used to find the moments of $X$ and the other constants in the formulas in Theorem 2.3. However, we shall leave these formulas to the reader and use a slightly different approach, where we use (4.16) for $-\varepsilon + \mathrm{i}t$ and let $\varepsilon \searrow 0$ to find $\psi(\mathrm{i}t)$ and its derivatives (at least for the $t$ that we need). Let

$$
(4.17) \qquad m(s) := \mathbb{E}\, Y_s = \sum_{k=1}^\infty (1 - e^{-s/2^k}), \qquad s > 0,
$$

and note that $0 \le m(s) \le s$. Further, since $Y_s$ is a sum of independent Bernoulli variables, $\operatorname{Var} Y_s \le \mathbb{E}\, Y_s$, and thus $\mathbb{E}\, Y_s(Y_s - 1) = \operatorname{Var} Y_s + m(s)^2 - m(s) \le m(s)^2 \le s^2$. Consequently, by a standard Taylor expansion of the probability generating function $D(s, u)$,

$$
D(s, u) = 1 + m(s)(u - 1) + O(s^2 |u - 1|^2), \qquad |u| \le 1.
$$

Evaluating (4.16) at $t = -\varepsilon + 2\pi n\mathrm{i}$ with $n \in \mathbb{Z}$ and $\varepsilon > 0$ thus yields

$$
\psi(t) = -\frac{t e^\varepsilon}{\ln 2} \int_0^\infty (1 + (e^{-\varepsilon} - 1) m(s) + O(s^2 \varepsilon^2)) s^{-t/\ln 2 - 1} e^{-s}\, ds,
$$

$$(4.18) \qquad = -\frac{t}{\ln 2}\left(e^{\varepsilon}\Gamma(-t/\ln 2) - \varepsilon \int_0^{\infty} m(s)s^{-t/\ln 2 - 1}e^{-s}\,ds\right) + O(\varepsilon^2)$$

$$= e^{\varepsilon}\Gamma(1 - t/\ln 2) + \varepsilon\frac{t}{\ln 2}g(-t/\ln 2) + O(\varepsilon^2),$$

where, for $\operatorname{Re} z > 0$, using (4.17),

$$g(z) := \int_0^{\infty} m(s)s^{z-1}e^{-s}\,ds$$

$$= \sum_{k=1}^{\infty}\int_0^{\infty}s^{z-1}(e^{-s} - e^{-(1+2^{-k})s})\,ds$$

$$= \sum_{k=1}^{\infty}(1 - (1 + 2^{-k})^{-z})\Gamma(z)$$

$$(4.19) \qquad = -\Gamma(z)\sum_{k=1}^{\infty}\sum_{j=1}^{\infty}\binom{-z}{j}2^{-jk}$$

$$= \sum_{j=1}^{\infty}(-1)^{j-1}\frac{\Gamma(z+j)}{j!\,(2^j - 1)};$$

we define $g(z)$ for $z \neq -1, -2, \dots$ by the latter sum. Letting $\varepsilon \to 0$ in (4.18) and recalling $t = -\varepsilon + 2\pi n\mathrm{i}$, we find

$$(4.20) \qquad \varphi(2\pi n) = \psi(2\pi n\mathrm{i}) = \Gamma(1 - 2\pi n\mathrm{i}/\ln 2)$$

as in several other examples, and (first taking $d/d\varepsilon$)

$$(4.21) \quad \psi'(2\pi n\mathrm{i}) = -\Gamma\left(1 - \frac{2\pi n}{\ln 2}\mathrm{i}\right) - \frac{1}{\ln 2}\Gamma'\left(1 - \frac{2\pi n}{\ln 2}\mathrm{i}\right) - \frac{2\pi n\mathrm{i}}{\ln 2}g\left(-\frac{2\pi n}{\ln 2}\mathrm{i}\right).$$

In particular, $\mathbb{E}\,X = \psi'(0) = \gamma/\ln 2 - 1 \approx -0.16725382272$. We similarly find

$$\mathbb{E}\,X^2 = \psi''(0) = \frac{\Gamma''(1)}{\ln^2 2} - 2\frac{\gamma}{\ln 2} + 1 - \frac{2g(0)}{\ln 2},$$

where

$$g(0) = \sum_{j=1}^{\infty}\frac{(-1)^{j-1}}{j(2^j - 1)} = \sum_{k=1}^{\infty}\ln(1 + 2^{-k}) \approx 0.86887665,$$

and thus

$$(4.22) \qquad \operatorname{Var} X = \frac{\pi^2}{6\ln^2 2} - \frac{2g(0)}{\ln 2} \approx 0.916666666667904.$$

As remarked in Example 2.11, $\operatorname{Var} X + \frac{1}{12}$ is extremely close to 1.

It follows from Theorem 2.1 and Remark 4.2 that

$$\mathbb{E}\,Y_s - \log_2 s = \mathbb{E}\,X_{\log_2 s} + o(1) = \gamma/\ln 2 - 1/2 + \beta_1(\log_2 s) + o(1)$$

(cf. [13], Answers 6.3-31 and 5.2.2-48, [9, 10]; recall also (4.17)).

For the variance, we find directly from the definition

$$\operatorname{Var} Y_s = \sum_{k=1}^{\infty} (1 - e^{-s/2^k}) e^{-s/2^k}$$

$$= \sum_{k=1}^{\infty} (e^{-s/2^k} - e^{-s/2^{k-1}}) = 1 - e^{-s},$$

since the sum telescopes. Thus, as $s \to \infty$, $\operatorname{Var} Y_s \to 1$. Hence, by Remark 4.2 (uniform square integrability is easily verified), $\operatorname{Var} X_\alpha = 1$ for every $\alpha$; see [9]. As discussed in Example 2.11, there is thus no oscillation in the variance.

For higher moments, see [9].

This example can also be generalized to bases other than 2 [2, 23]. Note, however, that the cancellation of the oscillations in the variance of $X_\alpha$ is special for base 2 [2] [and for $2^{1/k}$, $k \geq 1$, which does not make sense for Patricia tries but occurs for sampling from a geometric distribution $\operatorname{Ge}(p)$ with $q = 1 - p = 1/2^{1/k}$].

**5. Final remarks.** It is easily checked that if $X$ has the distribution function $F(x) := \mathbb{P}(X \leq x)$, then $X_\alpha$ has the distribution function

(5.1)                    $\mathbb{P}(X_\alpha \leq x) = F(\lfloor x + \alpha \rfloor - \alpha).$

Now suppose that (4.4) holds with some continuous function $F$, not necessarily a distribution function, and consider a subsequence such that $\{a_n\} \to \alpha \in [0, 1]$. Then,

$$\mathbb{P}(Y_n - a_n \leq x) = \mathbb{P}(Y_n \leq x + a_n)$$

$$= \mathbb{P}(Y_n \leq \lfloor x + a_n \rfloor)$$

$$= F(\lfloor x + a_n \rfloor - a_n) + o(1)$$

$$= F(\lfloor x + \{a_n\} \rfloor - \{a_n\}) + o(1)$$

$$\to F(\lfloor x + \alpha \rfloor - \alpha)$$

when $x + \alpha \notin \mathbb{Z}$. Hence $Y_n - a_n \xrightarrow{\mathrm{d}} X_\alpha$, as in (4.3), where $X_\alpha$ has the distribution given by (5.1). In particular, the right-hand side of (5.1) is a distribution function, even if $F$ is not (provided $\alpha$ is a limit point of $\{a_n\}$); this (if valid for all $\alpha$) is just equivalent to $F(x) \leq F(x + 1)$ for all $x$.

Note further that if $F$ has bounded variation, we can define $\varphi(t) := \int e^{\mathrm{i}tx} \, dF(x)$ and moments $\int x^m \, dF(x)$, and Theorems 2.1 and 2.3 still hold, with natural interpretations, by the same proofs. Thus, for example, results about asymptotics for moments of $Y_n$ can be obtained without knowing whether $F$ is a distribution function or not. (This is done, with different methods, in [17, 18, 25].)

Nevertheless, in each example where (4.4) holds for some $F$, we consider it to be an interesting question whether or not the limit $F$ is a distribution function. We have given several examples above where the answer is affirmative. One example where the problem is open is provided by *leader election* in the biased case ($p \neq 1/2$) [11]. Further, there are several examples where it seems that this question has not yet been studied; the interested reader can start by investigating the remaining examples in [17].

**Acknowledgments.** I would like to thank Dag Jonsson, Guy Louchard and Helmut Prodinger for helpful comments, and Robert Parviainen for finding an old calculation of mine.

DEPARTMENT OF MATHEMATICS
UPPSALA UNIVERSITY
PO BOX 480
S-751 06 UPPSALA
SWEDEN
E-MAIL: svante.janson@math.uu.se
URL: http://www.math.uu.se/~svante/